\documentclass[oneside,a4paper,11pt]{article}
\usepackage{amsfonts, amsbsy, amsmath, amsthm, amssymb, latexsym, verbatim, enumerate}
\usepackage{authblk}
\usepackage{graphicx,amssymb,amsmath}
\usepackage{amsthm}
\usepackage{amssymb}
\usepackage{latexsym}
\usepackage{longtable}
\usepackage{epsfig}
\usepackage{hhline}
\usepackage{hyperref}
\usepackage{aliascnt}
\newtheorem{theorem}{Theorem}[section]
\newaliascnt{lemma}{theorem}
\newtheorem{lemma}[lemma]{Lemma} 
\aliascntresetthe{lemma}
\newaliascnt{proposition}{theorem}
   
\aliascntresetthe{proposition}
\newaliascnt{corollary}{theorem}
\newtheorem{corollary}[corollary]{Corollary}
\aliascntresetthe{corollary}
\newaliascnt{question}{theorem}

\aliascntresetthe{question}
\newaliascnt{conjecture}{theorem}

\aliascntresetthe{conjecture}

\newtheorem*{theorem*}{Theorem}
\numberwithin{theorem}{section}
\newcommand{\pf}{{\it Proof:\quad}}

\newcommand{\dne}{\hfill $\Box$ \vspace{0.3cm}}

\title{Singular Graphs on which the Dihedral Group Acts Vertex Transitively }
\author{Ali Sltan Ali AL-Tarimshawy\\Department of Mathematices and Computer applications, College of Science, AL-Muthanna University\\ alisltan81@yahoo.com}

\begin{document}
\maketitle

\begin{abstract}
Let $\Gamma$ be a simple  connect graph on a finite vertex set $V$ and let $A$ be its adjacency matrix. Then $\Gamma$ is said to be \textit{singular} if and only if $0$ is an eigenvalue of $A.$ The \textit{nullity (singularity)} of $\Gamma,$ denoted by ${\rm null}(\Gamma),$ is the \textit{algebraic multiplicity} of the eigenvalue $0$ in the spectrum of $\Gamma.$   The general problem of characterising singular graphs is easy to state but it seems too difficult in this time. In this work, we investigate this problem for  finite graphs on which the  dihedral group $D_n$ acts vertex transitively as group of automorphisms. We determine the nullity of such graphs. We show that Cayley graphs over dihedral groups $D_{p^s} $  is non-singular if $|H \cap C_{p^s}|\neq |H \cap C_{p^s}b|$ and $|H|<p$ where $p$ is a prime number and $s \in \mathbb{N}.$
\end{abstract}
\section{\sc Introduction}

All graphs in this paper are  undirected, simple and finite. Let $\Gamma=(V,E)$ be a   graph on a vertex set $V$  together with a set $E$ of edges. Two vertices $v$ and $w$ are said  to be \textit{adjacent} in $\Gamma$ if and only if $\{v,w\}\in E.$ The \textit{adjacency matrix} of $\Gamma$ is the integer matrix with rows and columns indexed by the vertices of $\Gamma,$ such that the $a_{vw}$-entry is equal to $1$ if and only if $v \sim w$ and $0$ otherwise.  It  is denoted by $A.$  The \textit{characteristic polynomial} of $A$ is the polynomial $$C_A(x)=det(x I_n-A).$$ Here $I_n$ is the $(n\times n)$ identity matrix. It is clear that the roots of the characteristic polynomial of $A$ equal the eigenvalues of $A.$ Then $\Gamma$ is \textit{singular} if $A$ is singular. The \textit{spectrum }of $\Gamma$  consists of all eigenvalues  of $A$ and so $\Gamma$ is singular if and only if $0$ belongs to the spectrum of $\Gamma.$  The \textit{nullity} of  $\Gamma$ is the dimension of the null space of $\Gamma$ and  we denote this by ${\rm null}(\Gamma).$   Note $|V|={\rm null}(\Gamma) +r(\Gamma)$ where $r(\Gamma)$ is the rank  of $A.$   Hence  singular graphs have a non-trivial null space.

The problem of graph singularly first arose in structural chemistry in the context of H\"uckel Theory~\cite{graovac1972graph}  and the first mathematical paper on the subject appears to be  Collatz and  Sinogowitz~\cite{von1957spektren} in 1957 who asked for a classification of all finite non-singular graphs. Recently graph singularity has become relevant in other areas of mathematics too. This applies to the representation theory of finite groups in particular.  For instance, in \cite{muller2005some} M\"{u}ller and  Neunh\"{o}ffer compute the rank of a certain (extremely large) matrix that is connected to the Foulkes conjecture and the plethysm of symmetric groups. It is now evident that the problem of  Collatz and Singowitz can not be solved in its original version, it would simply not be feasible to determine all singular graphs.

Let $\mathbb{C}$ be the field of complex numbers. Let $G$ be a finite group. A subset $H$ of $G$ is called a \textit{connecting set} if  (i) $H^{-1}=\{h^{-1}: h \in H\}=H,$  (ii) $1_G \notin H$ and  (iii) $H$ generates $G.$ In this case we define a graph with vertex set $G$ and two vertices $u,v \in G$ are adjacent, $u \sim v$ if and only if $vu^{-1} \in H.$ This graph is called the \textit{Cayley graph} of $G$ with connecting set $H$ and it is denoted by $Cay(G,H).$  We denote the \textit{group algebra} of  $G$ over  $\mathbb{C}$ by $\mathbb{C}G$. Then  $\mathbb{C}G$  is the vector space over $\mathbb{C}$ with  basis $G$ and multiplication defined by extending the group multiplication linearly.

 A permutation $g$ of $V$ is an \textit{automorphism} of $\Gamma$ if the pair of vertices $(u^g,v^g)$ forms an edge in $\Gamma$ if and only if the pair of vertices $(u,v)$ forms an edge in $\Gamma.$ Here $u^g$ is the image of $u$ under the action of $g.$ The set of all automorphisms of $\Gamma$ forms a subgroup of the symmetric group on $V,$ called the \textit{automorphism group} of the graph $\Gamma.$ It is denoted by $Aut(\Gamma).$ In this note we consider the singularity problem for graphs which admit a group $G$ of automorphisms that acts transitively on the vertices of the graph. We specialize in particular to the case when $G$ is a dihedral group. In this situation the classification of singular graphs is complete.

To state our results as before  $\Gamma$ is  a graph with vertex set $V$ and let $G=D_{n}=\langle a,b\vert a^n=b^2=1,bab=a^{-1}\rangle\subseteq {\rm Aut}(\Gamma)$ be the dihedral group of order $2n.$ We suppose that $G$ acts transitively on $V.$ We denote the rotation group of $G$ by $C_{n}.$ Let $a$ be a generator of $C_n$ and let $b\in G\setminus C_n.$ 

We fix a vertex $v$ in $V$ and define the set 
$$H=\{g\in G :  \text{$v$ is adjacent to $v^g$\}}.$$ Write  $H=H' \cup H''$  where $H'= H\cap C_n$  and $H''=H\setminus H'.$ Next we define the polynomials

\begin{equation}\label{pap2_5}
\Psi'(x)=\sum\, x^{k}  \,\,\text{for $a^{k}\in H'$}
\end{equation}
 and
 \begin{equation}\label{pap2_6}
 \Psi''(x)=\sum\, x^{k}  \,\,\text{for $a^{k}b\in H''$}.
\end{equation}

\begin{lemma}\label{pap211} Let $G=D_n$ be the dihedral group of order $2n$ with cyclic subgroup $C_n$ of order $n.$ Suppose that  $G$ acts faithfully and transitively on some set  $V.$ Then  \\[5pt](i) $|V|=n$ and $C_n$ acts regularly on $V,$ or  \\[5pt] (ii)  $|V|=2n$ and $G$ acts regularly on $V$ while $C$ has two orbits on $V.$ 
\end{lemma}

\pf Let $F$ be the stabilizer of the point $v \in V$ and let $X:= F\cap C_n.$ Since $C_n$ is cyclic the only subgroup of $C_n$ of order  $|X|$ is $X$ itself, and in particular, $X$ is normal in $G.$ By a general lemma on permutation groups, $F^g=g^{-1} F g$ is the stabilizer of $v^{g}$ for any $g\in G,$ and so $X=X^g\subseteq F^g$ fixes $v^g$ for all $g.$ By transitivity, $X$ fixes any point of $V$ and  that means $X=1_G,$ as $F$ is faithful on $V.$ In particular, $|F|=1$ or $|F|=2.$  Since $|G|=|V||F|$ by the orbit stabilizer theorem we have the two options in the lemma. \dne

The first case,  if $|V|=n,$ then we have that  $C_n$ acts regularly on $V$ Hence   by Subidussi's  Theorem, see \cite[Lemma~3.7.2]{godsil2013algebraic} as a reference, we have that  $\Gamma$ is  a  Cayley graph over $C_n$ with a  connecting set $H'=\{g \in C_n : v \sim v^g\}.$ In this case we have that $A$ is a circulant matrix see the definition in Section~2. The second case, if $|V|=2n$ then $\Gamma$ is a Cayley graph over $G$ with connecting set \begin{equation}\label{pap2_4}
 H=\{ g \in G: v \sim v^g\} 
 \end{equation} and  hence  we have that $H=H'\cup H''.$ Thus  for this case $A$  has the following shape:
 
  \begin{equation}\label{pap2_3}
A= \bordermatrix{~ & C_n & C_nb \cr
                  C_n& M & N\cr
                  C_nb & N^t & M \cr} 
                  \end{equation} where $M,N$ are $n\times n$ matrices and $N^t$ is the transpose of $N.$

Note, if $G=D_n$ acts transitively on a set $V$ and suppose that  the kernel of this action $L$  is non-trivial. In this case  we have that $L \subseteq C_n$ and $G/L\cong D_{\tilde{n}}$ for some $\tilde{n}<n$ as  $L$ is a normal subgroup.

  In this paper, we  determine sufficient conditions for $\Gamma$ to be singular in both cases. We   prove that the eigenvalues of $A$ are those of $(M+N)$ together with those of $(M-N),$  when $\Gamma$ is a Cayley graph over dihedral groups see \autoref{pap26} and we  show that ${\rm null}(\Gamma)={\rm null}(M+N)+{\rm null}(M-N),$
see \autoref{pap28}. We  prove that $Cay(D_{p^{s}},H)$ is non-singular if the following conditions hold 
   
 (1) $|H\cap C_{p^s}|\neq |H\cap C_{p^s}b|;$
 
 (2) $|H|< p,$ 
 
 where $p$ is a prime number and $s \in \mathbb{N},$ see \autoref{pap210}.

\section{Circulant Matrices}
In this section we give the basic ideas and results for circulant matrices. We use the properties of such matrices in the next section to prove the main results of this work. 

Let $y=(y_0,y_1,...,y_{n-1})$ be a vector in $\mathbb{C}^n.$ An $n \times n$ \textit{circulant matrix } $Y$ takes the form $$Y = 
 \begin{pmatrix}
  y_0 & y_1 & \cdots & y_{n-1} \\
  y_{n-1} & y_0 & \cdots & y_{n-2} \\
  \vdots  & \vdots  & \ddots & \vdots  \\
 y_1 & y_2 & \cdots & y_0
 \end{pmatrix}$$ and it is denoted by $Y=circ(y).$ The elements of each row of $Y$ are identical to those of the previous row but are moved one position to the right and wrapped around. Note,  If $X,Y$ are  $n \times n$ circulant matrices, then $rX+sY$ is  $n\times n$ circulant matrix where $r,s$ are constant see \cite{daviscirculant} as a reference.
 
  The polynomial $$\Psi^*_Y(x)=y_0+y_1x+y_2x^2+...+y_{n-1}x^{n-1}$$ is called the \textit{associated polynomial} of matrix $Y,$ see \cite{daviscirculant} as a reference. Hence by the above polynomial we can determine the eigenvalues and the singularity of circulant matrices as shown in \autoref{pap23} and \autoref{pap24}.

      Let $\l$ be a positive  integer and let $\Omega_{\l}$ be the group of  $\l^{ th}$ \textit{roots of unity}, that is  $\Omega_{\l}=\{z\in \mathbb{C}\backslash\{0\} :z^{\l}=1\}.$ Then $\Omega_{\l}$ is a cyclic group of order $\l$ with  generator  $e^\frac{2\pi i}{\l}.$ Note this is not  the only generator of $\Omega_{\l},$  indeed any power $e^\frac{2\pi im}{\l}$ where $gcd(\l,m)=1$ is a generator too. A generator of $\Omega_{\l}$ is called \textit{a primitive $\l ^{th}$ root of unity}. Let $\Phi_{\l}(x)$ denote the ${\l}^{th}$ \textit{cyclotomic polynomial}. Then $\Phi_{\l}(x)$ is the unique irreducible integer polynomial with leading  coefficient $1$ so that $\Phi_{\l}(x)$ divides  $x^{\l}-1$ but  does not divide  of $x^k-1$ for any $k<{\l}$. Its roots are all primitive ${\l}^{th}$  roots of unity. So  $$\Phi_{\l}(x)=\prod_{1\leq m<{\l}}(x-e^\frac{2\pi im}{{\l}})$$ where $gcd(m,{\l})=1.$ \textit{Euler's totient function} of $\l$ is defined as the number of positive integer $\leq \l$ that is relatively prime to $\l$ and it is denoted by $\varphi(\l).$

 \begin{lemma}\label{pap23} \cite{lal2011non}
 Let $Y=circ(y)$ be an $n\times n$ circulant matrix. Then the eigenvalues of $Y$ equal $\Psi^*_Y(\omega^{i})$, for $i=0,1,...,n-1,$ where $\omega$ is a primitive $n^{th}$ root of unity. 
 \end{lemma}
   
 \begin{lemma}\label{pap24}\cite{kra2012circulant}
  Let $Y=circ(y)$ be an $n\times n$ circulant matrix. Then $Y$ is singular if and only if $\Phi_d(x)$ divides $\Psi^*_Y(x),$ for some divisor $d$ of $n$  where $\Phi_d(x)$ is a $d^{th} $ cyclotomic polynomial. Then the nullity of $Y$ is $\sum\varphi(d)$ where the sum is over all $d$ such that $\Phi_d(x)$ divides $\Psi^*_Y(x)$   and $\varphi(d)$ is the  Euler's totient function.
 \end{lemma}

\section{The Main Results of this work}
In this section we give sufficient conditions for $\Gamma=(V,E)$ to be singular where $V$ is a finite set and $G=D_n\subseteq Aut(\Gamma)$ acts faithfully and transitively on $V.$  We also determine the nullity of such graphs. We use the properties of the cyclotomic polynomial and circulant matrices to prove the main results in this paper.

Now we will discuss  some aspects about the Cayley graphs over dihedral groups which are related to our problem.  As before  $G=D_n=\langle a,b\vert a^n=b^2=1,bab=a^{-1}\rangle$ is the dihedral group of order $2n$ and  $C_n=\langle a\rangle$ is the cyclic group of order $n.$ Indicate a Cayley graph over $G$ and a connecting set $H$ by $Cay(G,H).$ Then  $A$ has the  shape is shown in \autoref{pap2_3}.

 Now we show  that $N=N^t.$  We define the following linear maps $$\beta:\mathbb{C}C_n\rightarrow \mathbb{C}C_nb,$$ it is given by $$\beta(a^{i})=\sum_{a^kb \in H''}ba^{-k}a^i$$ and $$\beta^t:\mathbb{C}C_nb\rightarrow \mathbb{C}C_n,$$ it is given by $$\beta^t(a^{\l}b)=\sum_{a^kb H''}ba^{-k}a^{\l}b.$$ Note $N$ and $ N^t$ represents $\beta$ and $\beta^t$ on the basis of $\mathbb{C}C_n$ and $\mathbb{C}C_nb$ respectively.

 Let $a^m,a^mb \in G.$ Then we have that

\begin{eqnarray}\label{pap2_1}
\beta(a^m)&=&\sum_{a^kb \in H''}ba^{-k}a^m\nonumber\\
&=&a^{k-m}b\nonumber\\
\end{eqnarray} and  
\begin{eqnarray}\label{pap2_2}
\beta^t(a^{m}b)&=&\sum_{a^kb \in H''}ba^{-k}a^{m}b\nonumber\\
&=& a^{k-m}.\nonumber\\
\end{eqnarray}

From this we conclude that $N=N^t.$ Hence $A$ has the following shape 
$$
A= \bordermatrix{~ & C_n & C_nb \cr
                  C_n& M& N\cr
                  C_nb & N & M\cr}.$$
                  
Now  we compute the eigenvalues  of $A$ as shown:
\begin{eqnarray}\label{pap2_7}
  C_A(\lambda)&=&{\rm det}(\lambda I_{2n}-A)\nonumber\\
  &=&{\rm det}
  \left( {\begin{array}{cc}
   \lambda I_n-M& N \\
   N &\lambda I_n- M \\
  \end{array} } \right) \nonumber\\
&=& {\rm det}\left( {\begin{array}{cc}
  \lambda I_n- M+N & \lambda I_n- M+N \\
   N& \lambda I_n- M\\
  \end{array} } \right)\nonumber\\
  &=&{\rm det}\left( {\begin{array}{cc}
   \lambda I_n- M+N & 0 \\
   N& \lambda I_n- (M+N) \\
  \end{array} } \right)\nonumber\\
  &=& {\rm det}(\lambda I_n-M+N)\cdot  {\rm det}(\lambda I_n- (M+N))\nonumber\\
 \end{eqnarray}
where $0$ is an $n\times n$ zero matrix. From this we deduce that  the eigenvalues of $A$ are those of $M+N$ and together with those of $M-N.$

\begin{theorem}\label{pap26}
Let $G=D_n=\langle a,b\vert a^n=b^2=1,bab=a^{-1}\rangle$ be the dihedral group of order $2n.$ Let $H$ be a connecting set for the graph $Cay(G,H).$ Then the eigenvalues of $A$ are those of $M+N$ together with $M-N.$
\end{theorem}

\pf It is clear by \autoref{pap2_7}.\dne

\begin{corollary} \label{pap28}
 Let $G=D_n=\langle a,b\vert a^n=b^2=1,bab=a^{-1}\rangle$ be the dihedral group of order $2n.$ Let $H$ be a connecting set for the graph $Cay(G,H).$ Then ${\rm null}(Cay(G,H))={\rm null}(M+N)+{\rm null}(M-N).$
\end{corollary}

 Clearly,  $M$ is a circulant matrix. Note  $\beta$ is also a circulant map as $\beta(a^i)=\sum_{a^kb \in H'' }a^{k-i}b$ from this we conclude that  if $V=(v_0,v_1,...,v_{n-1})$ is  the first row of $N,$ then the next  row of $N$ will move one position to the left, and wrapped around  and so on. From this we conclude that $M+N$ and $M-N$ are circulant matrices.
By the Cayley graph definition we have that all the elements in the connecting set $H$ of $Cay(G,H)$ are adjacent to $1_G.$ Hence  By \autoref{pap2_5} and \autoref{pap2_6} we have  the following polynomials  $$\Psi^*_M(x)=\Psi'(x)=\sum\, x^{k}  \,\,\text{for $a^{k}\in H'$}$$ and $$\Psi^*_N(x)=\Psi''(x)=\sum\, x^{k}  \,\,\text{for $a^{k}b\in H''$}$$ are associated to $M$ and $N$ respectively. Thus, the associated polynomial to $M+N$ and $M-N$ are $$\Delta^+(x)=\Psi'(x)+\Psi''(x)$$ and $$\Delta^-(x)=\Psi'(x)-\Psi''(x)$$ respectively.
 
 By this  we have the following result.
 
   \begin{theorem}\label{pap27}
 Let $G=D_n=\langle a,b\vert a^n=b^2=1_G,bab=a^{-1}\rangle$ be the dihedral group of order $2n.$ Let $H$ be a connecting set for the graph $Cay(G,H).$ Then $Cay(G,H)$ is singular if and only if either $\Phi_d(x)$ divides $\Delta^+(x)$ or $\Phi_d(x)$ divides $\Delta^-(x)$ for some divisors  $d$ of $n$ and ${\rm null}(Cay(G,H))=\sum \varphi(d)$ where the sum over all $d$ such that $\Phi_d(x)$ divides $\Delta^+(x)$ or $\Delta^-(x).$ 
 \end{theorem}

\begin{theorem} \label{pap210}
 Let $D_{p^s}=\langle a,b\vert a^{p^s}=b^2=1,bab=a^{-1}\rangle$ be the dihedral group of order $2p^s$ where $p$ is a prime number and $s \in \mathbb{N}.$ Let $H$ be a connecting set for the graph $Cay(D_{p^s},H).$ Then $Cay(D_{p^s},H)$ is non-singular if the following two conditions hold
 
 (1) $|H\cap C_{p^s}|\neq |H\cap C_{p^s}b|;$
 
 (2) $|H|< p.$
\end{theorem}
\pf (1) Suppose that the first condition is not held. In this case we have that the associated polynomial  to $M-N$ is $$\Delta^-(x)=\sum_{a^j \in H'}x^j-\sum_{a^{\l}b \in H''} x^{\l}.$$ Note we have that $\Delta^-(1)=0$ as $|H\cap C_{p^s}|= |H\cap C_{p^s}b|.$ From this we deduce that $\Phi_1(x)=x-1$ divides $\Delta^-(x)$ so by \autoref{pap24} we have that $M-N$ is singular with nullity $1.$ Hence by \autoref{pap28} we have that $Cay(D_{p^s},H)$ is singular with nullity $1.$ Otherwise if this condition is held then $\Gamma$ is non-singular.

(2) Now suppose that the first condition hold and $|H|<p.$ Suppose for contradiction $Cay(D_{p^s},H)$ is singular. In this case by \autoref{pap27} we have that either $\Phi_d(x)$ divides $\Delta^+(x)$ or $\Phi_d(x)$ divides $\Delta^-(x)$ for some divisor of $p^s.$ Note the divisors of $p^s$ are $1,p,p^2,...,p^s$ and by the first condition we have that $\Phi_1(x)$ is not divided $\Delta^+(x)$ or $\Delta^-(x).$ By this we have that at least one of the following Cyclotomic polynomials $\Phi_p(x),...,\Phi_{p^s}(x)$ divides $\Delta^+(x)$ or $\Delta^-(x).$   However   by the second condition we have that the number of terms of $\Delta^+(x)$ and $\Delta^-(x)$  is less than $p$ this give us a contradiction as the number of terms in each polynomial of $\Phi_p(x),...,\Phi_{p^s}(x)$ is $p^{th}$ terms. So $Cay(D_{p^s},H)$ is non-singular. \dne

Now we will comeback to our problem and the following result determine the singularity of $\Gamma.$

\begin{theorem} Let $\Gamma=(V,E)$ be an undirected graph without loops or multiple edges. Suppose that the dihedral group $G=D_{n}\subseteq {\rm Aut}(\Gamma)$ acts transitively and faithfully  on $V.$ Then \\[5pt]
(i)\,\, $|V|=n$ and $\Gamma$ is singular if and only if  $\Phi_{d}(x)$ divides $\Psi'(x)$ for some $d$ of $n$ and ${\rm null}(\Gamma)=\sum \varphi(d)$ where the sum is over all $d$ such that $\Phi_d(x)$ divides $\Psi'(x)$  or  \\[5pt]
(ii)\, $|V|=2n$ and $\Gamma$ is singular if and only if  $\Phi_{d}(x)$ divides $\Delta^+(x)$ or $\Delta^-(x)$ for some divisor $d$ of $n$ and ${\rm null}(\Gamma)=\sum \varphi(d)$ where the sum is over all $d$ such that $\Phi_d(x)$ divides $\Delta^+(x)$ or $\Delta^-(x).$  
\end{theorem}
\pf (i) \,\,Suppose that $|V|=n.$ In this case by Subidussi's Theorem  we have that $\Gamma$ is a Cayley graph over $C_n$ with connecting set $H'.$ Hence we have that $A$ is a circulant matrix  and  so by \autoref{pap24} we have that $A$ is singular if and only if $\Phi_d(x)$ divides $\Psi'(x)$ for some divisors $d$ of $n.$  

(ii)\,\, Suppose that $|V|=2n.$ Then by  by Subidussi's Theorem  we have that $\Gamma$ is a Cayley graph over $G$ with connecting set $H$ as  is shown in \autoref{pap2_4}. In this case we we have that $\Gamma=Cay(G,H)$ and the proof is completed   by applying \autoref{pap27}.


\end{document}